# Об энтропийно-подобных функционалах, возникающих в стохастической химической кинетике при концентрации инвариантной меры и в качестве функций Ляпунова динамики квазисредних

*Гасников А.В., Гасникова Е.В.*

В работе выявлен характер связи функции Ляпунова макросистемы, динамика которой определяется законами стохастической химической кинетики, с инвариантной мерой этой макросистемы, которая возникает на больших временах. Приведено необходимое и достаточное условие сводимости задачи поиска равновесия макросистемы (наиболее вероятного макросостояния инвариантной меры этой макросистемы) к задаче энтропийно-линейного программирования.

Библиография: 21 названий.

Предположим, что некоторая макросистема может находиться в различных состояниях, характеризуемых вектором $\vec{n}$ с неотрицательными целочисленными компонентами. Будем считать, что в системе происходят случайные превращения (химические реакции).

Пусть $\vec{n} \to \vec{n} - \vec{\alpha} + \vec{\beta}$, $(\vec{\alpha}, \vec{\beta}) \in J$ – все возможные типы реакций, где $\vec{\alpha}$ и $\vec{\beta}$ – вектора с неотрицательными целочисленными компонентами. Введем, следуя [1–3], интенсивность реакции:

$$\lambda_{(\vec{\alpha}, \vec{\beta})}(\vec{n}) = \lambda_{(\vec{\alpha}, \vec{\beta})}(\vec{n} \to \vec{n} - \vec{\alpha} + \vec{\beta}) = M^{1-\sum_i \alpha_i} K_{\vec{\beta}}^{\vec{\alpha}} \prod_{i: \alpha_i > 0} n_i \cdot \ldots \cdot (n_i - \alpha_i + 1),$$

где $K_{\vec{\beta}}^{\vec{\alpha}} \geq 0$ – константа реакции (поскольку допускается, что константа реакции может равняться нулю, то, не ограничивая общности, будем считать, что каждая реакция имеет обратную реакцию, возможно, с нулевой константой реакции); при этом часто считают $\sum_i n_i(t) \equiv M \gg 1$. Другими словами, $\lambda_{(\vec{\alpha}, \vec{\beta})}(\vec{n})$ – вероятность осуществления в единицу времени перехода $\vec{n} \to \vec{n} - \vec{\alpha} + \vec{\beta}$. На макроуровне это соответствует принципам химической кинетики (закону действующих масс Гульдберга–Вааге, 1864 [4, 5]), также говорят, что имеет место "приближение среднего поля" [1, 2]. Таким образом, динамика макросистемы задается линейной полугруппой (однородным дискретным марковским случайным процессом), инфинитезимальный оператор которой $\Lambda$ определяется интенсивностями реакций $\lambda_{(\vec{\alpha}, \vec{\beta})}(\vec{n})$: $d\vec{p}^T(t)/dt = \vec{p}^T(t)\Lambda$, где $[\Lambda]_{ii} = -\sum_{j:\, j \neq i}[\Lambda]_{ij}$. У вектора $\vec{p}(t)$ столько компонент, сколько существует допустимых конфигураций $\vec{n}$.

Везде в дальнейшем, не оговаривая этого, будем считать рассматриваемую марковскую динамику эргодической [6] (в смысле единственности стационарной меры).

Также отметим, что в отличие от стационарной меры (нормированной на единицу), инвариантная мера определяется с точностью до положительного множителя.

**Пример (модель П. и Т. Эренфестов в непрерывном времени, 1907 [7]).** На двух камнях сидят кузнечики. Всего кузнечиков $M \gg 1$. Каждый кузнечик независимо ни от чего в промежутке времени $[t, t+\Delta t)$, где $\Delta t$ – мало, а $t \geq 0$ – произвольно, перепрыгивает на другой камень с вероятностью $\lambda \Delta t + o(\Delta t)$, $\lambda > 0$. Введем вектор $\vec{n}(t) = (n_1(t), n_2(t))^T$, где $n_k(t)$ – число кузнечиков на $k$-м камне в момент времени $t \geq 0$. Тогда все возможные типы реакций задаются множеством: $J = \left\{ \left\{ \begin{pmatrix} 1 \\ 0 \end{pmatrix}, \begin{pmatrix} 0 \\ 1 \end{pmatrix} \right\}; \left\{ \begin{pmatrix} 0 \\ 1 \end{pmatrix}, \begin{pmatrix} 1 \\ 0 \end{pmatrix} \right\} \right\}$, при этом константы обеих "реакций" одинаковы и равны $\lambda$, а интенсивности реакций определяются формулами: $\lambda_{(1 \to 2)}(n_1, n_2) = \lambda n_1$, $\lambda_{(2 \to 1)}(n_1, n_2) = \lambda n_2$. Введем вектор $\vec{p}(t) = (p_0(t), ..., p_M(t))^T$, где $p_i(t)$ – вероятность того, что на первом камне в момент времени $t \geq 0$ ровно $i$ кузнечиков. Тогда инфинитезимальная матрица $\Lambda$ имеет следующий вид:

$$[\Lambda]_{ij} = \begin{cases} 0, & |i-j| > 1 \\ \lambda i, & j = i-1 \\ \lambda(M-i), & j = i+1 \\ -\lambda M, & j = i \end{cases}. \quad \Box$$

**Определение.** Условием унитарности или условием Штюкельберга–Батищевой–Пирогова (такое название было нам предложено пару лет назад В.В. Веденяпиным) называются следующие соотношения:

$$\exists \, \vec{\xi} > \vec{0}: \, \forall \, \vec{\beta} \to \sum_{\vec{\alpha}: (\vec{\alpha}, \vec{\beta}) \in J} K_{\vec{\beta}}^{\vec{\alpha}} \prod_j \xi_j^{\alpha_j} = \sum_{\vec{\alpha}: (\vec{\alpha}, \vec{\beta}) \in J} K_{\vec{\alpha}}^{\vec{\beta}} \prod_j \xi_j^{\beta_j}$$

или, что то же самое,

$$\exists \, \vec{\xi} > \vec{0}: \, \forall \, \vec{\alpha} \to \sum_{\vec{\beta}: (\vec{\alpha}, \vec{\beta}) \in J} K_{\vec{\beta}}^{\vec{\alpha}} \prod_j \xi_j^{\alpha_j} = \sum_{\vec{\beta}: (\vec{\alpha}, \vec{\beta}) \in J} K_{\vec{\alpha}}^{\vec{\beta}} \prod_j \xi_j^{\beta_j}. \qquad \text{(ШБП)}$$

Это условие обобщает хорошо известное в физике и экономике условие детального равновесия [7–10]:

$$\exists \, \vec{\xi} > \vec{0}: \, \forall \, (\vec{\alpha}, \vec{\beta}) \in J \to K_{\vec{\beta}}^{\vec{\alpha}} \prod_j \xi_j^{\alpha_j} = K_{\vec{\alpha}}^{\vec{\beta}} \prod_j \xi_j^{\beta_j}.$$

Так обычное условие детального равновесия будет выглядеть при дополнительном предположении, что инвариантная мера представима в виде прямого произведения распределений Пуассона: $\nu(\vec{n}) = \prod_i \xi_i^{n_i} e^{-\xi_i} / n_i!$.

Отметим, что в примере Эренфестов имеет место детальное равновесие с $\xi_1 = \xi_2 = 1$.

Следующая теорема [11, 12] уточняет результаты а) В.В. Веденяпина (2001) [4]; б) В.А. Малышева, С.А. Пирогова, А.Н. Рыбко (2004) [1, 2] на случай, когда не предполагается, что число состояний $\dim \vec{n}$ и число реакций $|J|$ не зависят от числа агентов $M$ (для примера Эренфестов – кузнечиков). Заметим, что есть довольно много примеров макросистем, когда это обобщение нужно [11, 12].

**Теорема 1. а)** $\langle \vec{\mu}, \vec{n}(t) \rangle \equiv \langle \vec{\mu}, \vec{n}(0) \rangle$ (inv) $\Leftrightarrow$ *вектор* $\vec{\mu}$ *ортогонален каждому вектору семейства* $\{\vec{\alpha} - \vec{\beta}\}_{(\vec{\alpha},\vec{\beta}) \in J}$. *Здесь* $\langle \cdot, \cdot \rangle$ – *обычно евклидово скалярное произведение.*

**б)** *Пусть выполняется условие (ШБП). Тогда*

1. *мера* $\nu(\vec{n}) = \prod_i \lambda_i^{n_i} e^{-\lambda_i} / n_i!$, *где* $\lambda_i = \xi_i^* M$, *а* $\vec{\xi}^*$ – *произвольное решение (ШБП), будет инвариантной относительно предложенной стохастической марковской динамики.*

2. *на множестве, определяемом уравнениями (inv), эта мера экспоненциально быстро концентрируется, с ростом* $M$, *в окрестности наиболее вероятного состояния, которое и принимается за положение равновесия макросистемы.*

3. *задача поиска наиболее вероятного макросостояния асимптотически эквивалентна задаче максимизации энтропийного функционала*:

$$\ln \nu(\vec{n}) \approx -\sum_i n_i \cdot \left( \ln(n_i / \lambda_i) - 1 \right) \tag{E}$$

*на множестве, задаваемом условием (inv).*

Решение задачи энтропийно-линейного программирования из ч. 3 п. б) теоремы 1 будем называть, следуя В.В. Веденяпину, экстремалью по Больцману.

Заметим, что в примере Эренфестов вектор $\vec{\mu} = (1,1)^T \perp \text{Lin} \{\vec{\alpha} - \vec{\beta}\}_{(\vec{\alpha},\vec{\beta}) \in J}$ определяет закон сохранения числа кузнечиков: $n_1(t) + n_2(t) \equiv M$, и это будет единственный закон сохранения. Стационарная (инвариантная) мера имеет вид:

$$\nu(n_1, n_2) = \nu(c_1 M, c_2 M) = M! \frac{(1/2)^{n_1}}{n_1!} \frac{(1/2)^{n_2}}{n_2!} = C_M^{n_1} 2^{-M} \simeq \frac{2^{-M}}{\sqrt{2\pi c_1 c_2}} \exp(-M \cdot H(c_1, c_2)),$$

где $H(c_1, c_2) = \sum_{i=1}^{2} c_i \ln c_i$. Кстати, сказать, из такого вида стационарной меры, можно получить, что если в начальный момент все кузнечики находились на одном камне, то математическое ожидание времени первого возвращения макросистемы в такое состояние будет порядка $2^M/M$ [7]. На примере этой модели можно говорить о том, что в макросистеме возврат к неравновесным макросостояниям вполне допустим, но происходить это может только через очень большое время (циклы Пуанкаре), так что нам может не хватить отведенного времени, чтобы это заметить (парадокс Цермело). Напомним, что описанный выше случайный процесс обратим во времени. Однако, как будет показано далее, наблюдается необратимая динамика относительной разности числа кузнечиков на камнях (парадокс Лошмидта).

Следующий пример, показывает, что условие (ШБП) является только достаточным условием инвариантности "пуассоновской" меры.

**Контрпример б) (С.А. Пирогов, 2004).** Рассмотрим систему уравнений химических реакций (константы реакций $K$ одинаковы):

$$A + B^* \xrightarrow{K} A^* + B, \quad B + C^* \xrightarrow{K} B^* + C, \quad C + A^* \xrightarrow{K} C^* + A, \text{ причем}$$

$$n_A(t) + n_{A^*}(t) \equiv n_B(t) + n_{B^*}(t) \equiv n_C(t) + n_{C^*}(t) \equiv N,$$

$$n_A(t) + n_B(t) + n_C(t) \equiv 3N/2.$$

Можно проверить, что "пуассоновская" мера

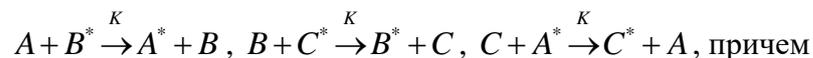

$$\nu(\vec{n}) \sim C_N^{n_A} \cdot C_N^{n_B} \cdot C_N^{n_C} \sim \underbrace{\left(1^{n_A} e^{-1}/n_A!\right) \cdot \ldots \cdot \left(1^{n_{C^*}} e^{-1}/n_{C^*}!\right)}_{\text{6 множителей}}$$

будет инвариантной, хотя условие (ШБП), очевидным образом, не выполняется. □

В виду только что рассмотренного примера, можно задаться вопросом: какое условие не только достаточно, но и необходимо для того, чтобы "пуассоновская" мера была инвариантной мерой макросистемы, и, как следствие, поиск равновесия сводился к поиску экстремали Больцмана в виде решения задачи энтропийно-линейного программирования? Оказывается, что условие (ШБП) и будет необходимым условием с некоторой оговоркой, которая будет приведена в следствие 2 к теореме 2.

Приведем другой подход, позволяющий получать во многом аналогичные результаты. Если в первом подходе мы сначала делали предельный переход по времени (стремится к бесконечности), а потом по числу агентов (стремится к бесконечности), то во втором подходе предлагается осуществить предельные переходы в обратном порядке.

Предположим теперь, что множество $J$ не зависит от $M$, и в начальный момент времени для любого $i$ существует предел

$$c_i(0) = \lim_{M \to \infty} n_i(0)/M. \qquad \text{(П)}$$

Тогда в произвольный момент времени $t > 0$ и для любого $i$ с вероятностью 1 существует предел (заметим, что $n_i(t)$ – случайные величины, тем не менее $c_i(t)$ – уже не случайные величины) $c_i(t) \stackrel{\text{п.н.}}{=} \lim_{M \to \infty} n_i(t)/M$. Описанный выше приём называется каноническим скейлингом. В результате такого скейлинга приходим к «динамике квазисредних» (терминология В. Вайдлиха [10]):

$$\frac{dc_i}{dt} = \sum_{(\vec{\alpha}, \vec{\beta}) \in J} (\beta_i - \alpha_i) K_{\vec{\beta}}^{\vec{\alpha}} \vec{c}^{\vec{\alpha}}, \quad \vec{c}^{\vec{\alpha}} = \prod_j c_j^{\alpha_j}. \qquad \text{(ДК)}$$

Это технически нетривиальное утверждение следует из результатов Троттера-Куртца (1986) [1, 13]. Отметим, что систему (ДК) можно получить и по-другому. А именно, как приближенную динамику средних $\bar{c}_i(t) = E[n_i(t)/M]$. Отметим также, что аффинное многообразие, задаваемое соотношениями (inv), с заменой $n$ на $c$, будет инвариантным многообразием системы (ДК) [4]. Для примера Эренфестов система (ДК) имеет вид:

$$\frac{dc_1}{dt} = \lambda(c_2 - c_1),$$

$$\frac{dc_2}{dt} = \lambda(c_1 - c_2).$$

Можно показать, что если выполняются условия (ШБП), то любая траектория (ДК) сходится к неподвижной точке. Какой именно, зависит, вообще говоря, от «точки старта»; но можно сказать и точнее: к той единственной неподвижной точке из семейства неподвижных точек, которая принадлежит аффинному многообразию (inv), инвариантному относительно (ДК).[1] Для этого, вводится (минус) энтропия (E): $H(\vec{c}) = \sum_i c_i \cdot (\ln(c_i/\xi_i) - 1)$ и показывается, что она является функцией Ляпунова для системы (ДК). Этот результат, по-видимому, впервые был получен Я.Г. Батищевой и В.В. Веденяпиным (2001) [4].

---

[1] Стоит заметить, что аттрактор системы (ДК), по-видимому, в общем случае может быть сколь угодно сложным множеством [2, 10].

Для примера Эренфестов: $H(\vec{c}) = \sum_i c_i \cdot \ln(c_i)$, $c_1(t) + c_2(t) \equiv 1$ (inv). Глобально устойчивое положение равновесия (экстремаль Больцмана) определяется из условия:

$$\vec{c}^* = \begin{pmatrix} 1/2 \\ 1/2 \end{pmatrix} = \arg \min_{\substack{c_1+c_2=1 \\ \vec{c} \geq 0}} H(\vec{c}).$$

Обратим внимание, что инвариантная мера (при каноническом скейлинге) "породила" функцию Ляпунова (см. п. б) ч. 3 теоремы 1). Это не случайно. Подобные закономерности наблюдаются для рассматриваемых моделей и без условия (ШБП), и даже без предположения о том, что инвариантная мера концентрируется около единственного положения равновесия, то есть аттрактором может быть множество куда как более сложной структуры.

**Теорема 2.** *Пусть выполняется условие (П) и инвариантная мера представляется в виде:*

$$\nu(\vec{n}) = M \exp\left(-M \cdot \left(H(\vec{n}/M) + o(1)\right)\right), \vec{n} \in (\text{inv}), \quad M \gg 1. \quad (\text{ИМ})$$

*Тогда если $H(\vec{c})$ строго выпуклая функция, то $H(\vec{c})$ – функция Ляпунова системы (ДК).*

**Доказательство.** Введем производящую функцию:

$$F(t, \vec{s}) = \sum_{\vec{n}} P(\vec{n}(t) = \vec{n}) \cdot \vec{s}^{\vec{n}}, \quad |\vec{s}| \leq \vec{1}, \quad \vec{s}^{\vec{n}} = s_1^{n_1} \cdot s_2^{n_2} \cdot \ldots,$$

на которую можно выписать следующее уравнение в частных производных [14]:

$$\frac{\partial F(t, \vec{s})}{\partial t} = \sum_{(\vec{\alpha}, \vec{\beta}) \in J} \left(\vec{s}^{\vec{\beta}} - \vec{s}^{\vec{\alpha}}\right) K_{\vec{\beta}}^{\vec{\alpha}} M^{1-\sum_i \alpha_i} \frac{\partial^{\alpha_1+\alpha_2+\ldots} F(t, \vec{s})}{\partial s_1^{\alpha_1} \cdot \partial s_2^{\alpha_2} \cdot \ldots}. \quad (\text{УЧП})$$

Отметим, что если взять частные производные $\partial/\partial s_i$ от обеих частей (УЧП), полагая $\vec{s} = (1, \ldots, 1)^T$, то получим в пределе при $M \to \infty$ систему (ДК).

Поскольку по формуле Эйлера–Маклорена [15]

$$F(\infty, \vec{s}) \approx M \int e^{M\left(\langle \overrightarrow{\ln s}, \vec{\zeta} \rangle - H(\vec{\zeta})\right)} d\vec{\zeta},$$

то [16]

$$M^{-(\alpha_1+\alpha_2+\ldots)} \frac{\partial^{\alpha_1+\alpha_2+\ldots} F(\infty, \vec{s})}{\partial s_1^{\alpha_1} \cdot \partial s_2^{\alpha_2} \cdot \ldots} \approx \frac{\vec{\zeta}(\vec{s})^{\vec{\alpha}}}{\vec{s}^{\vec{\alpha}}} C(\vec{s}, M),$$

где $\vec{\zeta}(\vec{s})$ определяется, притом единственным образом, из системы $\overrightarrow{\ln s} = \operatorname{grad} H(\vec{\zeta})$, а $C(\cdot) \neq 0$ – не зависит от $\vec{\alpha}$ (воспользовались методом Лапласа асимптотического оценивания интеграла, зависящего от параметров, и его производных по этим параметрам). Следовательно,

$$0 \equiv \sum_{(\vec{\alpha},\vec{\beta})\in J} \left(\vec{s}^{\vec{\beta}} - \vec{s}^{\vec{\alpha}}\right) K_{\vec{\beta}}^{\vec{\alpha}} \frac{\vec{\zeta}(\vec{s})^{\vec{\alpha}}}{\vec{s}^{\vec{\alpha}}} = \sum_{(\vec{\alpha},\vec{\beta})\in J} \left(e^{\langle \vec{\beta}-\vec{\alpha},\,\text{grad}\,H(\vec{\zeta}(\vec{s}))\rangle} - 1\right) K_{\vec{\beta}}^{\vec{\alpha}} \vec{\zeta}(\vec{s})^{\vec{\alpha}} \geq$$

$$\geq \sum_{(\vec{\alpha},\vec{\beta})\in J} \langle \vec{\beta}-\vec{\alpha},\,\text{grad}\,H(\vec{\zeta}(\vec{s}))\rangle K_{\vec{\beta}}^{\vec{\alpha}} \vec{\zeta}(\vec{s})^{\vec{\alpha}} = \left.\frac{dH(\vec{c})}{dt}\right|_{\vec{c}=\vec{\zeta}(\vec{s})}$$

– полная производная функции $H(\vec{c})$ в силу системы (ДК) в точке $\vec{\zeta}(\vec{s})$. □

**Следствие 1.** *Пусть выполняются условия (П) и (ИМ) с гладкой функцией $H(\vec{c})$ (не обязательно выпуклой). Тогда аттрактор системы (ДК) принадлежит множеству $\text{Arg}\min_{\vec{c}\in(\text{inv})} H(\vec{c})$, то есть при больших значениях времени произвольная траектория системы (ДК) достаточно близка к этому множеству.*

**Следствие 2.** *Пусть выполняются условия (П) и (ИМ) со строго выпуклой функцией $H(\vec{c})$. Тогда для того, чтобы в условии (ИМ) $H(\vec{c}) = \sum_i c_i \cdot (\ln(c_i/\xi_i) - 1)$ для почти все $\vec{c}(0)$ необходимо, чтобы выполнялось условие (ШБП).*

Для доказательства следствия 2 заметим, что $\overrightarrow{\ln s} = \text{grad}\,H(\vec{\zeta}) \Rightarrow s_i = \zeta_i/\xi_i$ и

$$0 \equiv \sum_{(\vec{\alpha},\vec{\beta})\in J} \left(\vec{s}^{\vec{\beta}} - \vec{s}^{\vec{\alpha}}\right) K_{\vec{\beta}}^{\vec{\alpha}} \frac{\vec{\zeta}(\vec{s})^{\vec{\alpha}}}{\vec{s}^{\vec{\alpha}}} \stackrel{s_i=\zeta_i/\xi_i}{=} \sum_{(\vec{\alpha},\vec{\beta})\in J} \vec{\zeta}^{\vec{\beta}} \frac{K_{\vec{\beta}}^{\vec{\alpha}} \vec{\xi}^{\vec{\alpha}}}{\vec{\xi}^{\vec{\beta}}} - \sum_{(\vec{\alpha},\vec{\beta})\in J} \vec{\zeta}^{\vec{\alpha}} K_{\vec{\beta}}^{\vec{\alpha}} =$$

$$= \sum_{\vec{\beta}} \vec{\zeta}^{\vec{\beta}} \left[ \sum_{\vec{\alpha}:(\vec{\alpha},\vec{\beta})\in J} \frac{K_{\vec{\beta}}^{\vec{\alpha}} \vec{\xi}^{\vec{\alpha}}}{\vec{\xi}^{\vec{\beta}}} - \sum_{\vec{\alpha}:(\vec{\alpha},\vec{\beta})\in J} K_{\vec{\alpha}}^{\vec{\beta}} \right] = \sum_{\vec{\alpha}} \vec{\zeta}^{\vec{\alpha}} \left[ \sum_{\vec{\beta}:(\vec{\alpha},\vec{\beta})\in J} \frac{K_{\vec{\alpha}}^{\vec{\beta}} \vec{\xi}^{\vec{\beta}}}{\vec{\xi}^{\vec{\alpha}}} - \sum_{\vec{\beta}:(\vec{\alpha},\vec{\beta})\in J} K_{\vec{\beta}}^{\vec{\alpha}} \right] \equiv 0.$$

Поскольку последние два соотношения должны выполняться для любых $\vec{\zeta} \in (\text{inv})$, то для почти всех $\vec{c}(0)$, на основе которых определяется множество (inv), должно выполняться условие (ШБП). □

**Следствие 3.** *Пусть выполняется условие (П)*

$$\exists\, c_i(0) > 0:\; n_i(0) = c_i(0)M + o(M),\; M \gg 1$$

*и (ИМ) со строго выпуклой функцией $H(\vec{c})$. Тогда время выхода макросистемы на равновесие есть $T = \text{O}(\ln M)$.*

Мы не приводим здесь доказательства следствия 3 (см. [17]). Упомянем лишь, что в основе доказательства лежат современные варианты изопериметрического неравенства Чигера [18].

В заключение хочется заметить, что в общем случае рассматриваемые в работе предельные переходы не перестановочны, и "жизнь" системы (ДК) не всегда всецело определяется только линейными законами сохранения, унаследованными ею при скейлинге из стохастической динамики. В результате скейлинга могут возникать новые (нелинейные) законы сохранения (первые интегралы). Все это хорошо иллюстрирует следующая макросистема.

**Пример (модель хищник–жертва, 1925 [4, 5, 19]).** Обозначим число зайцев ("жертв") в лесу $[Ж]$, а число волков ("хищников") $[X]$. Считаем, что волки и зайцы бродят случайно по лесу в поисках пищи. Зайцы плодятся согласно закону Мальтуса, но рост популяции зайцев лимитируют волки. Пищей волков являются зайцы (потенциально ограниченный ресурс), пища зайцев – трава (не ограниченный ресурс). Если волк "встречает" зайца, то он его съедает и размножается, то есть рождается новый волк. Если волк долгое время не встречает зайца – он погибает от голода. Зайцы выбывают только за счет встречи с волками. Формализуем сказанное в введенной ранее терминологии:

$$Ж \xrightarrow{K_1} 2Ж, \quad Ж + X \xrightarrow{K_2} 2X, \quad X \xrightarrow{K_3} 0,$$

$$\lambda_1(Ж, X) = K_1 \cdot [Ж], \quad \lambda_2(Ж, X) = K_2 \cdot [Ж][X], \quad \lambda_3(Ж, X) = K_3 \cdot [X],$$

$$\begin{aligned}\frac{dc_Ж}{dt} &= K_1 c_Ж - K_2 c_Ж c_X, \\ \frac{dc_X}{dt} &= K_2 c_Ж c_X - K_3 c_X.\end{aligned} \quad \text{(ДК)}$$

Система (ДК) имеет (нелинейный) первый интеграл:

$$K_3 \ln c_Ж + K_1 \ln c_X - K_2 (c_Ж + c_X) \equiv \text{const},$$

возникшей в результате скейлинга, то есть не присущий исходной стохастической марковской динамике. Таким образом, типичными траекториями системы (ДК) будут замкнутые циклы, отвечающие колебанию численности хищников и жертв. Система (ДК) хорошо известна в математической биологии и называется системой Лотки–Вольтера.

Однако стохастическая динамика имеет поглощающие состояния: без хищников (зайцев бесконечно много), и зайцев и волков нет. Существуют не нулевые вероятности (в сумме меньшие единицы) "свалиться" в одно из этих состояний (чем больше система, тем эти вероятности меньше). □

Обратим внимание на популярную статью [19], в которой рассматривается океан, представленный двумерной клеточной решеткой, по клеткам "случайно" блуждают акулы и скумбрии, … Исследование на больших промежутках времени при термодинамическом предельном переходе (размер океана стремится к бесконечности, пропорционально с численностями акул и скумбрий) такой динамики приводит к приведенным выше выводам (формулам). Таким образом, введенные в начале статьи формулы для интенсивностей реакций получают вполне естественное "биологическое" обоснование.

Отметим также, что в биологических приложениях, а также в популяционной теории игр с помощью используемого в статье формализма получается содержательная интерпретация мальтузианских параметров (функций). Таким образом, возникает завязка рассматриваемых в статье результатов на популяционную теорию игр, эволюционную оптимальность, дарвиновский отбор [20, 21].

## СПИСОК ЦИТИРОВАННОЙ ЛИТЕРАТУРЫ